\documentclass{article}

\usepackage{todonotes}
\usepackage{fullpage}
\usepackage{amsmath,amssymb,amsthm}
\usepackage[pagebackref=true]{hyperref}
\usepackage{tikz-cd}
\usepackage{dynkin-diagrams}
\usepackage{xtab}
\allowdisplaybreaks

\numberwithin{equation}{section}

\theoremstyle{definition}
\newtheorem{definition}{Definition}[section]
\newtheorem{theorem}[definition]{Theorem}
\newtheorem{proposition}[definition]{Proposition}

\newtheorem{remark}[definition]{Remark}

\newtheorem{problem}[definition]{Problem}

\DeclareMathOperator{\chern}{c}
\DeclareMathOperator{\ch}{ch}
\DeclareMathOperator{\td}{td}
\DeclareMathOperator{\Sym}{Sym}

\DeclareMathOperator{\rank}{rank}

\DeclareMathOperator{\Hom}{Hom}

\DeclareMathOperator{\res}{res}
\newcommand{\pt}{\mathrm{pt}}

\newcommand{\Hline}[1]{\noalign{\hrule height #1}}

\title{Elliptic genera of complete intersection Calabi--Yau 17-folds in $F_4$-Grassmannians}
\author{Kobayashi Kenta}
\date{}

\begin{document}
\maketitle

\begin{abstract}
  Motivated by the derived invariance problem
  of elliptic genera,
  we construct a non-birational pair of Calabi--Yau complete intersection 17-folds
  in $F_4$-Grassmannians
  with distinct Chern numbers
  and identical elliptic genera.
\end{abstract}

\section{Introduction}

A Calabi--Yau manifold is a compact K\"{a}hler manifold such that $c_1(M) \in \operatorname{H}_*(M; \mathbb{R})$ equals $0$.
Calabi--Yau manifolds are of great interest in mathematics and mathematical physics.
In particular,
Calabi--Yau $3$-folds play an important role in various fields of mathematics derived from superstring theory.
For 2 dimensions and below,
the classification have been made,
but for 3 dimensions and above,
it is not known whether the number of deformation equivalence classes is finite or not.

Elliptic genera were originally defined by Ochanine and Landweber--Stong
for oriented manifolds
\cite{MR895567,MR948178}.
An elliptic genus
in the sense of Ochanine
is defined as a genus that vanishes for the projective bundle of an even-dimensional complex vector bundle on an oriented closed manifold,
and is characterized as a genus
whose logarithm takes the form
\begin{equation}
  g(x) = \int_0^x \frac{dt}{\sqrt{1 - 2 \delta t^2 + \epsilon t^4}}
\end{equation}
where $\delta$, $\epsilon$ are constants.
By specializing $\delta$ and $\epsilon$,
we obtain the signature ($\delta = \epsilon = 1$)
and the $\hat{A}$-genus ($\delta = -1/8, \epsilon = 0$).
For an elliptic genus $\varphi$, we have
\begin{equation}
  \delta = \varphi(\mathbb{C}P^2), \quad \epsilon = \varphi(\mathbb{H}P^2)
\end{equation}
where $\mathbb{H}$ is the quaternion field.

By substituting the following formal power series in $q$
for $\delta$ and $\epsilon$,
we obtain the universal elliptic genus defined by Landweber--Stong \cite{MR970290}:
\begin{align}
  \delta_A
   & =
  - \frac{1}{8}
  - 3 \sum_{n \geq 1} \left( \sum_{d|n, \, 2 \nmid d } d \right) q^n, \\
  \epsilon_A
   & =
  - \sum_{n \geq 1} \left( \sum_{d|n, \, 2 \nmid (n/d)} d^3 \right) q^n.
\end{align}
The signature and $\hat{A}$-genus are the values of this genus at the cusp.
Furthermore,
this is a modular form for the group
\begin{equation}
  \Gamma_0(2) \coloneqq \left\{ \begin{pmatrix}
    a & b \\
    c & d
  \end{pmatrix}
  \in \operatorname{SL}_2(\mathbb{Z})
  \, \middle| \,
  c \equiv 0 \mod 2
  \right\}.
\end{equation}

Hirzebruch \cite{MR981372} and Witten \cite{MR970278} defined an elliptic genus for almost complex manifolds with its first Chern class divisible by $N$,
which is a modular form of the group
\begin{equation}
  \Gamma_1(N) \coloneqq \left\{
  \begin{pmatrix}
    a & b \\
    c & d
  \end{pmatrix}
  \in \operatorname{SL}_2(\mathbb{Z})
  \, \middle| \,
  \begin{pmatrix}
    a & b \\
    c & d
  \end{pmatrix}
  \equiv
  \begin{pmatrix}
    1 & b \\
    0 & 1
  \end{pmatrix}
  \mod N
  \right\}.
\end{equation}
Krichever \cite{MR1048541} and H\"{o}hn \cite{arxiv:math/0405232}
introduced an elliptic genus which can be specialized into the Hirzebruch--Witten elliptic genus.
The Krichever--H\"{o}hn elliptic genus provides weak Jacobi forms for $\operatorname{SU}$-manifolds.

By the rigidity theorem of Bott--Taubes \cite{MR954493},
the Ochanine--Landweber--Stong elliptic genus $\varphi$
has the following strict multiplicativity property:
For a fiber bundle $E$
which has a compact connected Lie group as the structure group,
an oriented closed manifold $B$ as the base space,
and a closed spin manifold $F$ as the fiber,
we have
\begin{equation}
  \varphi(E) = \varphi(F)\varphi(B).
\end{equation}
The Krichever--H\"{o}hn elliptic genus
has a similar strict multiplicativity property
for a fiber bundle which has a compact connected Lie group as the structure group,
a closed stable complex manifold as the base space,
and a closed $\operatorname{SU}$-manifold as the fiber \cite{MR1048541, arxiv:math/0405232}.
This property plays a crucial role in the proof of the theorem of Totaro
to be mentioned later.

In this paper,
we deal with a variant of the Krichever--H\"{o}hn elliptic genus
originating from physics \cite{MR985505,MR1264785,MR1463039}.
The precise definition is given in Section \ref{sc:ell}.
Hereafter,
we refer to this variant simply as the elliptic genus.

A pair of $\mathbb{Q}$-Gorenstein projective algebraic varieties $X$ and $X'$ are said to be \emph{K-equivalent}
if there exist a smooth projective algebraic variety $Y$
and birational maps $f \colon Y \rightarrow X$ and $f' \colon Y \rightarrow X'$
such that the pullbacks of the canonical divisors, $f^* K_X$ and ${f'}^* K_{X'}$, are $\mathbb{Q}$-linearly equivalent.
By Totaro \cite{MR1765709},
the elliptic genus provides a ring isomorphism
from the quotient of the cobordism ring $\mathit{MSU}_* \otimes \mathbb{Q}$
by the ideal generated by differences of varieties connected by a classical flop
to the ring of weak Jacobi forms over $\mathbb{Q}$.
Moreover,
this statement holds even if we replace the classical flops with K-equivalence \cite{MR1953295,MR1949645}.
In other words,
the elliptic genus is universal
among characteristic numbers
which are invariant under K-equivalences.

With the DK hypothesis (see e.g.~\cite{MR3838122} and its references)
in mind,
the following problem arises naturally (cf.~\cite[Section 7]{MR2180406}):

\begin{problem} \label{problem:dtoe}
Is the elliptic genus invariant under derived equivalences?
\end{problem}

In this paper,
we construct a pair of Calabi--Yau manifolds
that are
not K-equivalent
but are expected to be derived equivalent,
and confirm the equality of their elliptic genera.
If the dimension of the manifold is less than $12$
or equal to $13$,
then the elliptic genus is determined by the $\chi_y$-genus (and therefore by the Hodge numbers)
\cite[Theorem 2.7]{MR1757003}.
Since Hodge numbers are expected to be invariant under derived equivalences,
a counter-example to Problem \ref{problem:dtoe}, if any,
would be in high dimensions.

In dimension three,
Calabi--Yau complete intersections in exceptional Grassmannians have been exhaustively studied by Ito--Miura--Okawa--Ueda and Benedetti \cite{arxiv:1606.04076,bened2018}.
We construct $17$-dimensional Calabi--Yau complete intersections
in $F_4$-Grassmannians similar to the $G_2$-cases in their works.
Let
$\dynkin[label]{F}{4}$
be the Dynkin diagram of type $F_4$
and $\omega_i$ be the fundamental weight corresponding to the $i$-th vertex.
The main theorem of this paper is as follows:

\begin{theorem}\label{th:main}
  Let $G$ be the connected and simply connected simple complex Lie group of type $F_4$,
  and $X_1 \coloneqq G/P_1$ (resp. $X_2 \coloneqq G/P_2$) be the homogeneous space corresponding to the crossed Dynkin diagram $\dynkin{F}{*x**}$ (resp. $\dynkin{F}{**x*}$).
  Let further $\mathcal{E}_1$ (resp. $\mathcal{E}_2$) be the $G$-equivariant vector bundle on $X_1$ (resp. $X_2$) associated with the weight $\omega_2+\omega_3$,
  and $Y_1$ (resp. $Y_2$) be the zero locus of a general section of $\mathcal{E}_1$ (resp. $\mathcal{E}_2$).
  Then the following hold:
  \begin{enumerate}
    \item \label{it:CY} $Y_1$ and $Y_2$ are $17$-dimensional Calabi--Yau manifolds.
    \item \label{it:chern} $Y_1$ and $Y_2$ have different Chern numbers. In particular, they are not diffeomorphic.
    \item \label{it:non-birat} $Y_1$ and $Y_2$ are not birationally equivalent.
    \item \label{it:ell} $Y_1$ and $Y_2$ have the same elliptic genus.
  \end{enumerate}
\end{theorem}

In Section \ref{sc:CYCI},
we discuss the construction of the Calabi--Yau complete intersections $Y_1$ and $Y_2$,
and prove Theorems \ref{th:main}.\ref{it:CY},
\ref{th:main}.\ref{it:chern},
and \ref{th:main}.\ref{it:non-birat}.
In Section \ref{sc:ell},
we clarify the definition of the elliptic genus adopted in this paper
and prove \ref{th:main}.\ref{it:ell}.
The proofs of these theorems depend
on a SageMath package
\cite{arXiv:2306.06660}
developed for this purpose.
In addition,
the proof of Theorem \ref{th:main}.\ref{it:non-birat} relies on a Mathematica package by Makoto Miura \cite{arxiv:1606.04076}.

\subsection*{Acknowledgements}

I would like to express my deep gratitude to my supervisor, Kazushi Ueda, for raising this issue, providing numerous advice and carefully checking this manuscript.
The Calabi-Yau manifolds discussed in this paper were constructed using Makoto Miura's Mathematica package,
which is also used to prove that their Picard numbers are 1.
I am very grateful to Makoto Miura and his package.
I also appreciate Yuichi Enoki for telling me Equations \eqref{eq:phi0,1} and \eqref{eq:phi-2,1}.
This work is supported by Leading Graduate Course for Frontiers of Mathematical Sciences and Physics.

\section{$F_4$-Grassmannians and their complete intersections}\label{sc:CYCI}

Let $\mathfrak{g}$ be a complex semi-simple Lie algebra,
and
$G$ be the corresponding connected and simply connected Lie group.
Fix a Cartan subgroup $H$ of $G$
with the Lie algebra
$\mathfrak{h}$.
For $\alpha \in \mathfrak{h}^\vee$,
we set
\begin{equation}
  \mathfrak{g}_\alpha \coloneqq \left\{
  v \in \mathfrak{g} \,\middle|\,
  \forall h \in \mathfrak{h}, \ [h, v] = \alpha(h)v
  \right\},
\end{equation}
so that we have
\begin{equation}
  \mathfrak{g} = \mathfrak{h} \oplus \bigoplus_{\alpha \in \Delta} \mathfrak{g}_{\alpha}
\end{equation}
where
$
  \Delta \coloneqq
  \left\{ \alpha \in \mathfrak{h}^\vee \mid \mathfrak{g}_\alpha \neq \left\{ 0 \right\}, \  \alpha \neq 0\right\}
$
is the set of roots.
We choose a set
$
  \mathcal{S} \coloneqq \left\{ \alpha_1, \ldots, \alpha_r \right\}
  \subset \Delta
$
of \emph{simple roots},
so that the set $\Delta$ of roots is decomposed into the disjoint union
of the set $\Delta^+$ of positive roots
and the set $\Delta^-$ of negative roots.
The Weyl group $W_G$ is
generated by reflections through the hyperplanes orthogonal to simple roots.
The Dynkin diagram is a graph
whose vertices correspond to simple roots,
and edges (with decorations) carry information
of the inner products of simple roots.

A subgroup $P$ of $G$ is said to be parabolic
if $G/P$ is projective.
The set of conjugacy classes of parabolic subgroups
is bijective with
the set of subsets of simple roots of $G$.
For a subset
$
  \mathcal{S}_{\mathfrak{p}}
  \subset
  \mathcal{S}
$,
the Lie algebra $\mathfrak{p}$ of the corresponding parabolic subgroup is given by
\begin{align}
  \mathfrak{p} & \coloneqq
  \mathfrak{l}
  \oplus
  \mathfrak{n},            \\
  \mathfrak{l} & \coloneqq
  \mathfrak{h}
  \oplus
  \bigoplus_{\alpha \in (\operatorname{span} \mathcal{S}_{\mathfrak{p}}) \cap \Delta}
  \mathfrak{g}_\alpha,     \\
  \mathfrak{n} & \coloneqq
  \bigoplus_{\alpha \in \Delta^+ \setminus (\operatorname{span} \mathcal{S}_{\mathfrak{p}})}
  \mathfrak{g}_\alpha
\end{align}
where
$\mathfrak{l}$
is called the \emph{Levi part}
of
$\mathfrak{p}$,
and
$\mathfrak{n}$
is called the \emph{nilpotent part}
of $\mathfrak{p}$.
Conjugacy classes of parabolic subgroups
can be represented by crossed Dynkin diagrams:
We cross out the vertices
corresponding to roots
which are not included in $\mathcal{S}_{\mathfrak{p}}$.

For each simple root $\alpha_i$,
let
$\omega_i$
be the corresponding fundamental weight.
The weight lattice $\Lambda$
of $\mathfrak{g}$ is freely generated
by the set of fundamental weights.

A weight
$
  \lambda = \sum_{i=1}^{\rank G} \lambda_i \omega_i \in \Lambda \otimes \mathbb{R}
$
is \emph{integral}
if $\lambda \in \Lambda$.
It is
\emph{$\mathfrak{p}$-dominant}
if $\lambda_i \geq 0$ for all $i$ such that $\alpha_i \in \mathcal{S}_{\mathfrak{p}}$.
The highest weight representation of the Levi part $\mathfrak{l}$
corresponding to a $\mathfrak{p}$-dominant weight $\lambda$
lifts to a representation of $P$
if and only if $\lambda$ is integral
(see e.g.~\cite[Section 3.1]{MR1038279}).
We write the irreducible representation of $P$
associated with $\lambda$ as $V_{\lambda}^P$.

There exists an equivalence
\begin{equation}
  \begin{split}
    \operatorname{Rep}_{\mathbb{C}}(P) & \rightarrow \operatorname{Vect}_G(G/P)\\
    V & \mapsto \mathcal{E}_V^P \coloneqq G \times_P V
  \end{split}
\end{equation}
from the category of representations of a parabolic subgroup $P$
to the category of $G$-equivariant vector bundles over $G/P$.
The equivariant vector bundle associated with
an irreducible representation $V_{\lambda}^P$
is denoted by
\begin{equation}
  \mathcal{E}_{\lambda}^P \coloneqq \mathcal{E}_{V_{\lambda}^P}^{\vee}.
\end{equation}
This vector bundle is globally generated if and only if $\lambda$ is $\mathfrak{g}$-dominant.

Let $B$ be a Borel subgroup of $G$
and $G_0$ be the compact real form of $G$,
so that
$G/B \cong G_0/T$
where $T \coloneqq B \cap G_0$ is a maximal torus of $G_0$.
The canonical isomorphism
$T \cong \Hom(\Lambda, U(1))$
induces
isomorphisms
\begin{equation}
  H^*_G(G/B)
  \cong H^*_{G_0}(G_0/T)
  \cong H^*(BT)
  \cong \Sym \Lambda,
\end{equation}
where $BT$ is the classifying space of $T$.
The first Chern class $c_1(\mathcal{L}_{\omega_i})$
of the line bundle
$
  \mathcal{L}_{\omega_i} \coloneqq \mathcal{E}^{P}_{\omega_i}
$
on $G/B$
corresponds to
$\omega_i \in \Sym \Lambda$
under this isomorphism.

Similarly,
we have $G/P \cong G_0/P_0$
where
$P_0 \coloneqq G_0 \cap P$,
and the rational $G$-equivariant cohomology of $G/P$ can be calculated as
\begin{equation}
  H^*_G(G/P; \mathbb{Q})
  \cong H^*_{G_0}(G_0/P_0; \mathbb{Q})
  \cong H^*(BP_0; \mathbb{Q})
  \cong (\Sym \Lambda_{\mathbb{Q}})^{W_{P}}
\end{equation}
where $BP_0$ is the classifying space of $P_0$ and
$W_{P} \cong W_{P_0}$ is the Weyl group of $P$.

The tangent bundle $T_{G/P}$ of $G/P$ is associated
with the representation of $P$ on $\mathfrak{g}/\mathfrak{p}$ by the adjoint action:
\begin{equation}
  T_{G/P} \cong \mathcal{E}_{\mathfrak{g}/\mathfrak{p}}.
\end{equation}
Hence
the $G$-equivariant total Chern class of $G/P$ is given by
\begin{equation}
  c_G(G/P) = c_G(\mathcal{E}_{\mathfrak{g}/\mathfrak{p}}) = \prod_{\mu \in \Delta(\mathfrak{g}/\mathfrak{p}^{\vee})}(1 + \mu).
\end{equation}
Note that the set $\Delta(\mathfrak{g}/\mathfrak{p})$
of weights of $\mathfrak{g}/\mathfrak{p}$
is equal to $\Delta^- \setminus \operatorname{span} \mathcal{S}_{\mathfrak{p}}$.

Now assume that $G$ of type $F_4$ and
let $P_1$ and $P_2$ be the parabolic subgroups
corresponding to the crossed Dynkin diagrams
$\dynkin[label]{F}{*x**}$ and $\dynkin[label]{F}{**x*}$
respectively.
For each $i \in \{1, 2\}$ and
the $\mathfrak{g}$-dominant weight $\lambda = \omega_2 + \omega_3$,
the equivariant vector bundle $\mathcal{E}_i \coloneqq \mathcal{E}^{P_i}_\lambda$
on the homogeneous space $X_i \coloneqq G/P_i$ is globally generated.
The zero $Y_i$ of a general section of $\mathcal{E}_i$
is a smooth complete intersections
by a generalization of Bertini's theorem
\cite[Theorem 1.10]{MR1201388}.

Since $Y_i$ is a complete intersection,
one has an exact sequence
\begin{equation}
  0
  \rightarrow T_{Y_i}
  \rightarrow \iota^* T_{X_i}
  \rightarrow \iota^* \mathcal{E}_i
  \rightarrow 0,
\end{equation}
where
$\iota \colon Y_i \hookrightarrow X_i$
is the inclusion.
Due to the additivity and naturality of Chern classes, as well as the fact that $\iota^* \colon H^*(X_i) \rightarrow H^*(Y_i)$ is a ring homomorphism,
we have
\begin{equation}
  c(T_{Y_i})
  = \frac{c(\iota^* T_{X_i})}{c(\iota^* \mathcal{E}_i)}
  = \iota^* \left( \frac{c(T_{X_i})}{c(\mathcal{E}_i)} \right).
\end{equation}
Here,
$c(\iota^* \mathcal{E}_i)$ has an inverse in $H^*(Y_i)$.
In general,
for any vector bundle $\mathcal{E}$ on $Y_i$,
\begin{equation}
  \prod_{i = 1}^d (1 - x_i + x_i^2 - \cdots + (- x_i)^d)
\end{equation}
is the inverse element of $c(\mathcal{E}) = \prod_{i = 1}^d(1 + x_i)$.

Since the Euler class of the normal bundle of
$\iota \colon Y_i \hookrightarrow X_i$
is
$\iota^*(\prod_{\nu \in \Delta(V_{\lambda}^{\vee})} \nu)$,
the composite
$\iota_! \circ \iota^*$
is multiplication by $\prod_{\nu \in \Delta(V_{\lambda}^{\vee})} \nu$.
Define a homogeneous element $C_p(Y_i) \in \Sym \Lambda$ of degree $p$ by
\begin{equation}
  \sum_{p = 0}^\infty C_p(Y_i)
  = \frac{c_G(T_{X_i})}{c_G(\mathcal{E}_i)}
  = \frac{\prod_{\mu \in \Delta(\mathfrak{g}/\mathfrak{p}^{\vee})}(1 + \mu)}
  {\prod_{\nu \in \Delta(V_{\lambda}^{\vee})}(1 + \nu)}
  \in
  \widehat{\Sym \Lambda}.
\end{equation}
Then one has
\begin{equation}
  \int_{Y_i} c_{i(1)}\cdots c_{i(n)}
  =
  \int_{X_i} \left( C_{i(1)}(Y_i) \cdots C_{i(n)}(Y_i) \prod_{\nu \in \Delta(V_{\lambda}^{\vee})}\nu \right).
\end{equation}

\begin{proof}[{Proof of Theorem \ref{th:main}.\ref{it:CY}}]
  One has $\dim G = 52$ and
  $\dim P_i = 32$,
  so that
  $\dim X_i = 20$.
  Since $\rank V_{\lambda}^{P_i} = 3$,
  one has $\dim Y_i = 17$.
  Furthermore,
  using SageMath,
  we find that $C_1(Y_i)$,
  which is the degree-one part of
  $\frac{c_G(T_{X_i})}{c_G(\mathcal{E}_i)}$, vanishes
  (see \texttt{ancillary.ipynb}
  available as an ancillary file to the arXiv version of this paper).
  By taking the non-equivariant limit
  and pulling back,
  the first Chern class of $Y_i$ also vanishes.
\end{proof}

Let
$
  f \colon F_i \hookrightarrow X_i
$
be the inclusion
of the fixed locus
$
  F_i \coloneqq X_i^{T}
$.
We have the commutative diagram
\begin{equation}
  \begin{tikzcd}[row sep=15mm,column sep=20mm]
    H^{*}_T(X) & H^{*}_T(\mathrm{pt}) \\
    H^{*}_T(F_i) & H^{*}_T(F_i)
    \arrow["\int_{F_i}"',from=2-2, to=1-2]
    \arrow["f_!"',shift left=0,from=2-2, to=1-1]
    \arrow["f^{*}",from=1-1, to=2-1]
    \arrow["e_T(N_{F_i/X_i}) \cup "',hook',from=2-2, to=2-1]
    \arrow["\int_X",from=1-1, to=1-2]
  \end{tikzcd}
\end{equation}
where $e_T(N_{F_i/X_i})$ is the equivariant Euler class
of the normal bundle of $F_i$ in $X_i$.
There exists a $W_G$-equivariant bijection
from $W_G/W_{P_i}$
to
$
  F_i
$
sending
$
  [w] \in W_G/W_{P_i}
$
to
$
  wP_i \in F_i
$,
so that
$H^*_T(F_i)$ can be identified with
$\oplus_{[w] \in W_G/W_P} H^*_T(\pt)$,
and
$\int_{F_i}$ is the summation
$
  \oplus_{[w] \in W_G/W_P} H^*_T(\pt)
  \to
  H^*_T(\pt).
$
For any
$
  x \in H_G^*(X_i; \mathbb{Q})
  \cong (\Sym \Lambda_{\mathbb{Q}})^{W_P}
$,
we have
\begin{equation}
  f^* \circ \res (x) = \bigoplus_{[w] \in W_G / W_P} w \cdot x
\end{equation}
where
$
  \res \colon H_G^*(X_i; \mathbb{Q}) \to H_T^*(X_i; \mathbb{Q})
$
is the restriction
and
$w \cdot x$
is the action of
$w \in W_G$
on
$x \in \Sym \Lambda_{\mathbb{Q}}$
induced from the action of $W_G$
on $\Lambda$.

Since the fixed locus $F_i$ is discrete,
the normal bundle $N_{F_i/X_i}$
at $w P_i \in F_i$ is given by
\begin{equation}
  T_{w P_i} X_i \cong (\mathfrak{g}/\mathfrak{p})^w.
\end{equation}
The superscript on the right-hand side
shows that the $T$-action is twisted
by conjugation by $w$.
Hence
the Euler class at $w P_i$ is given by
\begin{equation}
  e_T \left(N_{F_i/X_i,wP_i} \right)
  = \prod_{\alpha \in \Delta(\mathfrak{g}/\mathfrak{p})} w \cdot \alpha,
\end{equation}
and
we obtain the localization formula
\begin{equation}\label{eq:loc}
  \int_{X_i} x
  =
  \sum_{
    [w] \in W_G / W_P
  }
  \frac{w \cdot x}{\prod_{\alpha \in \Delta(\mathfrak{g}/\mathfrak{p})} w \cdot \alpha}
\end{equation}
in the localized ring $H_T^*(\pt)_{e(N_{F_i/X_i})}$,
on which the SageMath module \cite{arXiv:2306.06660} is based.

\begin{proof}[{Proof of Theorem \ref{th:main}.\ref{it:chern}}]
  This is an immediate consequence of the table
  of Chern numbers of $Y_1$ and $Y_2$
  given in Appendix \ref{app:chernnum}
  obtained by using SageMath.
\end{proof}

\begin{proof}[{Proof of Theorem\ref{th:main}.\ref{it:non-birat}}]
  By using the Mathematica package developed by Makoto Miura
  (see \texttt{hodge.nb}
  available as an ancillary file to the arXiv version of this paper),
  we have
  \begin{equation}
    \operatorname{Pic}(Y_1)_{\mathbb{Q}} \cong
    \operatorname{Pic}(Y_2)_{\mathbb{Q}} \cong
    \mathbb{Q}.
  \end{equation}
  The integrals of the 17-th power of the Chern classes of the ample line bundles $L_i = \wedge^3 \mathcal{E}_i$ are
  $
    2^{13} \cdot 3^2 \cdot 5^{18} \cdot 13 \cdot 17 \cdot 3413
  $ and
  $
    2^7 \cdot 3^2 \cdot 5 \cdot 7^{18} \cdot 13^3 \cdot 17
  $,
  respectively (as seen in \texttt{ancillary.ipynb}).
  Since
  they cannot be related by a multiple of the $17$-th power of an integer,
  the 17-th power of the ample generators of $Y_i$ are distinct.
  Since birational Calabi--Yau manifolds are connected by a sequence of flops \cite{MR2426353},
  and
  a flop induces
  an isomorphism of Picard groups
  preserving the intersection numbers,
  the Calabi--Yau manifolds $Y_1$ and $Y_2$ are not birational.
\end{proof}

\begin{remark}
  The complete intersections of equivariant vector bundles corresponding to the weights $\omega_1+\omega_2$ on the Grassmannian variety of type $G_2$ corresponding to
  $\dynkin[label]{G}{x*}$ and $\dynkin[label]{G}{*x}$
  are Calabi--Yau 3-folds
  which are not birational \cite{arxiv:1606.04076}
  but are derived equivalent \cite{MR3830796}.
  Moreover,
  the derived equivalence of these complete intersection Calabi--Yau $3$-folds follows
  from the derived equivalence of the total spaces of the dual bundles
  of the corresponding equivariant vector bundles \cite[Section 5]{MR3959275}.
  Similarly,
  the derived equivalence of $Y_1$ and $Y_2$ follows from the derived equivalence of the total spaces of $\mathcal{E}_1^\vee$ and $\mathcal{E}_2^\vee$.
  Since the total spaces of $\mathcal{E}_1^\vee$ and $\mathcal{E}_2^\vee$ are related by a flop,
  the derived equivalence of $Y_1$ and $Y_2$ follows from the DK hypothesis.
  In addition,
  the argument of the proof of \cite[Lemma 2.1]{MR4110724}
  gives the relation
  \begin{equation}
    ([Y_1] - [Y_2]) \cdot \mathbb{L}^2 = 0
  \end{equation}
  in the Grothendieck ring of varieties,
  where
  $
    \mathbb{L} \coloneqq [\mathbb{A}^1]
  $
  is the class of the affine line.
\end{remark}

\section{Elliptic Genus} \label{sc:ell}

A \emph{complex ($\mathbb{Q}$-)genus} is
a ring homomorphism $\varphi \colon MU_* \rightarrow R$ from the complex cobordism ring $MU_*$ to a $\mathbb{Q}$-algebra $R$.
Complex genera are in one-to-one correspondence
with formal power series $Q(x)$ over $R$
satisfying $Q(0) = 1$
called their \emph{characteristic series}
(see \cite{MR1189136}).
Let
$
  c\left( M \right) = \prod_{i=1}^{\dim M}\left( 1 + x_i \right)
$
be the Chern class of a closed stable complex manifold $M$.
For a formal power series $Q(x) = 1 + a_1 x + a_2 x^2 + \cdots \in R[\![x]\!]$,
let
\begin{equation}
  \prod_{i = 1}^{\dim M} Q(x_i)
  =
  \sum_{r} K_r\left( c_1(M), c_2(M), \ldots c_r (M) \right)
\end{equation}
be the decomposition
into the sum of homogeneous parts $K_r$ of degree $r$.
Then the corresponding genus is given by
\begin{equation}
  \varphi_Q\left( M \right)
  \coloneqq
  K_{\dim M} \left( c_1(M), c_2(M), \ldots c_{\dim M} (M) \right)[M]
  \in R.
\end{equation}
Conversely,
given a genus $\varphi$,
we define its \emph{logarithm} as
\begin{equation}
  g(x) = \sum_{n = 0}^\infty
  \frac{\varphi(\mathbb{P}^n)}{n + 1}x^{n + 1}
  \in R[\![x]\!].
\end{equation}
Let $f(x)$ be the inverse series of $g(x)$
(in the sense that $g(f(x))=x$).
Then
the characteristic series is obtained as
$Q_\varphi(x) \coloneqq x/f(x)$.

The \emph{elliptic genus} is the complex genus
whose characteristic series is
\begin{align}
  Q(x)
   & =
  \frac{x\theta_1(\tau, \frac{x}{2\pi i} - z)}{\theta_1(\tau, \frac{x}{2\pi i})} \\
   & =
  y^{-1/2}
  \frac{x}{1 - e^{-x}}
  \prod_{n=1}^{\infty}
  \frac{
  (1 - yq^{n - 1}e^{-x})(1 - y^{-1}q^ne^{x})
  }{
  (1 - q^ne^{-x})(1 - q^ne^x)
  }
\end{align}
where
$
  R = \mathbb{Q}(y^{1/2})[\![q]\!]
$,
$q = e^{2\pi i \tau}$,
$y = e^{2 \pi i z}$, and
\begin{equation}
  \theta_1(z, \tau)
  =
  - q^{1/8}
  (2\sinh \pi z)
  \prod_{l = 1}^{\infty}(1 - q^l)
  \prod_{l=1}^{\infty}(1 - q^l y)(1 - q^l y^{-1}).
\end{equation}
For a $d$-dimensional complex manifold $M$,
we have
\begin{align}
  \chi(M, q, y)
   & =
  \int_M Q(x_1) \cdots Q(x_{d}) \\
   & =
  y^{-d/2}\int_M \td(M) \ch \left( \bigotimes_{n \geq 1}
  \left(
    \wedge_{-yq^{n-1}} T^*_M
    \otimes \wedge_{-yq^n T_M}
    \otimes
    \Sym_{q^n} T_M^*
    \otimes
    \Sym_{q^n}T_M\right)
  \right)
\end{align}
where
\begin{align}
  \Sym_q V
   & \coloneqq
  \bigoplus_{n = 0} ^{\infty}
  \left(
  \Sym^n V
  \right)
  \cdot
  q^n,         \\
  \wedge_q V
   & \coloneqq
  \bigoplus_{n = 0}^{\infty}
  \left(
  \wedge^n V
  \right)
  \cdot
  q^n.
\end{align}

A \emph{weak Jacobi form} of weight $k \in \mathbb{Z}$ and index $m \in \frac{1}{2}\mathbb{Z}$ is a function $\tilde{\phi} \left( \tau, z \right)$ with the Fourier series expansion
\begin{equation}
  \tilde{\phi} \left( \tau, z \right)
  =
  \sum_{n \geq 0}
  \sum_{r}
  c\left( n, r \right)
  q^ny^r
  =
  \sum_{n \geq 0}
  \sum_{r}
  c\left( n, r \right)
  \mathrm{e}^{2 \pi \mathrm{i} n \tau}
  \mathrm{e}^{2 \pi \mathrm{i} r z}
\end{equation}
which satisfies the following equations:
\begin{align}
  \tilde{\phi}
  \left(
  \frac{a \tau + b}{c \tau + d},
  \frac{z}         {c \tau + d}
  \right)
   & =
  (c \tau + d)^k
  \mathrm{e}^{\frac{2 \pi \mathrm{i} cz^2}{c \tau + d}}
  \tilde{\phi}
  \left( \tau, z \right),
   & \begin{pmatrix}
       a & b \\
       c & d
     \end{pmatrix}
   & \in \mathrm{SL}_2{\mathbb{Z}}, \\
  \tilde{\phi}
  \left(
  \tau,
  z + \lambda \tau + \mu
  \right)
   & =
  (-1)^{2m(\lambda + \mu)}
  \mathrm{e}
  ^{-2 \pi \mathrm{i} m(\lambda^2 \tau + 2 \lambda z)}
  \tilde{\phi}
  \left( \tau, z \right),
   & \left( \lambda, \mu \right)
   & \in \mathbb{Z}^2.
\end{align}

The space of weak Jacobi forms of even weight and integral index can be described as follows:

\begin{proposition}[{\cite[Theorem 9.3]{MR781735}}]
  The space
  $\tilde{J}_{\mathrm{even}, *} \coloneqq \bigoplus_{k \in 2\mathbb{Z}, \, m \in \mathbb{Z}}\tilde{\mathrm{J}}_{k, m}$
  of weak Jacobi forms of even weight and integral index
  is isomorphic to
  $\mathrm{M}_*\left[\tilde{\phi}_{-2, 1}, \tilde{\phi}_{0, 1}\right]$,
  where
  $\mathrm{M}_*$
  is the ring of modular forms
  (which is the polynomial ring
  generated by the normalized Eisenstein series
  of weight $4$ and $6$),
  \begin{align}\label{eq:phi0,1}
    \tilde{\phi}_{0, 1}(q, y)
    = 4 \left\{\left(\frac{\theta_2(q, y)}{\theta_2(q, 1)}\right)^2
    + \left(\frac{\theta_3(q, y)}{\theta_3(q, 1)}\right)^2
    + \left(\frac{\theta_4(q, y)}{\theta_4(q, 1)}\right)^2\right\}
  \end{align}
  is a weak Jacobi form of weight $0$ and index $1$, and
  \begin{align}\label{eq:phi-2,1}
    \tilde{\phi}_{-2, 1}(q, y)
     & = - \frac{\theta_1(q, y)^2}{\eta(q)^6}
  \end{align}
  is a weak Jacobi form of
  weight $-2$ and index $1$.
\end{proposition}

The space of weak Jacobi forms
of weight $0$ and half-integral index
is described as follows:

\begin{proposition}[{\cite[Proof of Theorem 2.6]{MR1757003}}]
  Let $d$ be an odd integer.
  Multiplication by
  \begin{equation}
    \tilde{\phi}_{0, \frac{3}{2}}(q, y)
    \coloneqq \frac{\theta_1(q, y^2)}{\theta_1(q, y)}
  \end{equation}
  gives an isomorphism
  $
    \tilde{J}_{0, \frac{(d-3)}{2}} \rightarrow \tilde{J}_{0, \frac{d}{2}}
  $
  from the space of weak Jacobi forms of weight $0$ and index $d/2$ to the space of weak Jacobi forms of weight $0$ and index $(d-3)/2$.
\end{proposition}

The elliptic genus of a Calabi--Yau manifold
is a weak Jacobi form:

\begin{proposition}[{\cite[Theorem 2.2]{MR1757003}}]
  The elliptic genus of a Calabi--Yau $d$-folds
  is a weak Jacobi form of index $d/2$ and weight $0$.
\end{proposition}

Now we prove Theorem \ref{th:main}.\ref{it:ell}:

\begin{proof}[{Proof of Theorem \ref{th:main}.\ref{it:ell}}]
  The space of weak Jacobi forms of weight $0$ and index $17/2$ is $8$-dimensional,
  and its basis is given by
  \begin{equation} \label{eq:basis}
    \left\{
    \begin{aligned}
      v_1 & = \tilde{\phi}_{0, \frac{3}{2}} \tilde{\phi}_{0, 1}^7,                                                  &
      v_2 & = \tilde{\phi}_{0, \frac{3}{2}} \mathrm{E}_4 \tilde{\phi}_{0, 1}^5 \tilde{\phi}_{-2, 1}^2,              &
      v_3 & = \tilde{\phi}_{0, \frac{3}{2}} \mathrm{E}_6 \tilde{\phi}_{0, 1}^4 \tilde{\phi}_{-2, 1}^3,              &
      v_4 & = \tilde{\phi}_{0, \frac{3}{2}} \mathrm{E}_4^2\tilde{\phi}_{0, 1}^3 \tilde{\phi}_{-2, 1}^4,               \\
      v_5 & = \tilde{\phi}_{0, \frac{3}{2}} \mathrm{E}_4 \mathrm{E}_6 \tilde{\phi}_{0, 1}^2 \tilde{\phi}_{-2, 1}^5, &
      v_6 & = \tilde{\phi}_{0, \frac{3}{2}} \mathrm{E}_6^2 \tilde{\phi}_{0, 1} \tilde{\phi}_{-2, 1}^6,              &
      v_7 & = \tilde{\phi}_{0, \frac{3}{2}} \mathrm{E}_4^3 \tilde{\phi}_{0, 1} \tilde{\phi}_{-2, 1}^6,              &
      v_8 & = \tilde{\phi}_{0, \frac{3}{2}} \mathrm{E}_4^2 \mathrm{E}_6 \tilde{\phi}_{-2, 1}^7
    \end{aligned}
    \right\}.
  \end{equation}
  A calculation using SageMath shows that
  their coefficient vectors up to first order in $q$
  are linearly independent
  (see \texttt{ancillary.ipynb}),
  so a weak Jacobi form of weight $0$ and index $17/2$
  is determined by its $q$-expansion up to first order.
  Using SageMath again,
  we compute the coefficients up to the first order of $q$ of the elliptic genus
  of a 17-dimensional Calabi--Yau in terms of the Chern numbers,
  and express the elliptic genus in terms of the basis \eqref{eq:basis}.
  The result is shown in Appendix \ref{app:ell17}.
  By combining this with the table of Chern numbers in Appendix \ref{app:chernnum},
  one finds that
  the elliptic genera of both $Y_1$ and $Y_2$ are given
  using this basis as
  \begin{multline}
    - \frac{523623513480701}{2985984} v_1
    + \frac{193611909253757}{331776} v_2
    - \frac{348708636989665}{1492992} v_3
    - \frac{442170341015857}{995328} v_4\\
    + \frac{22462744370909}{82944} v_5
    - \frac{5975835708317}{186624} v_6
    + \frac{17023636880707}{331776} v_7
    - \frac{10255680346099}{497664} v_8,
  \end{multline}
  and Theorem \ref{th:main}.\ref{it:ell} is proved.
\end{proof}

\appendix

\section{The elliptic genus of a Calabi--Yau 17-fold} \label{app:ell17}

The elliptic genus of a Calabi--Yau 17-fold
is expressed using Chern numbers as follows:
\begin{align*}
   & \chi\left(M^{17}, q, y\right)                                                                                                                                                                                                                          \\
   & = \int_M \frac{1}{71663616} c_{17} \cdot v_1                                                                                                                                                                                                           \\
   & - \int_M \left(  \frac{1}{119439360} c_{3} c_{14} + \frac{1}{119439360} c_{2} c_{15} - \frac{1}{19906560} c_{17}  \right)
  \cdot v_2                                                                                                                                                                                                                                                 \\
   & - \int_M \bigg( \frac{1}{250822656} c_{2} c_{3} c_{12} + \frac{1}{250822656} c_{2}^{2} c_{13} + \frac{1}{250822656} c_{5} c_{12} - \frac{1}{125411328} c_{4} c_{13}                                                                                    \\
   & + \frac{1}{125411328} c_{3} c_{14} - \frac{1}{250822656} c_{2} c_{15} + \frac{5}{250822656} c_{17} \bigg) \cdot v_3                                                                                                                                    \\
   & - \int_M \bigg( \frac{1}{597196800} c_{2}^{2} c_{3} c_{10} + \frac{1}{597196800} c_{2}^{3} c_{11}
  + \frac{1}{597196800} c_{3} c_{4} c_{10} + \frac{1}{597196800} c_{2} c_{5} c_{10}
  - \frac{1}{199065600} c_{2} c_{4} c_{11}                                                                                                                                                                                                                  \\
   & + \frac{1}{298598400} c_{2} c_{3} c_{12} - \frac{1}{597196800} c_{2}^{2} c_{13} - \frac{1}{597196800} c_{7} c_{10} + \frac{1}{199065600} c_{6} c_{11} - \frac{1}{119439360} c_{5} c_{12}                                                               \\
   & + \frac{1}{119439360} c_{4} c_{13} + \frac{1}{66355200} c_{3} c_{14} - \frac{11}{597196800} c_{2} c_{15} + \frac{13}{298598400} c_{17} \bigg) \cdot v_4                                                                                                \\
   & -\int_M \bigg( \frac{1}{1532805120} c_{2}^{3} c_{3} c_{8} + \frac{1}{1532805120} c_{2}^{4} c_{9} + \frac{1}{4598415360} c_{3}^{3} c_{8} + \frac{1}{766402560} c_{2} c_{3} c_{4} c_{8} + \frac{1}{1532805120} c_{2}^{2} c_{5} c_{8}                     \\
   & - \frac{1}{4598415360} c_{2} c_{3}^{2} c_{9} - \frac{1}{383201280} c_{2}^{2} c_{4} c_{9} + \frac{1}{766402560} c_{2}^{2} c_{3} c_{10} - \frac{1}{1532805120} c_{2}^{3} c_{11} - \frac{1}{1532805120} c_{4} c_{5} c_{8}                                 \\
   & - \frac{1}{1532805120} c_{3} c_{6} c_{8} - \frac{1}{1532805120} c_{2} c_{7} c_{8} + \frac{1}{766402560} c_{4}^{2} c_{9} + \frac{1}{4598415360} c_{3} c_{5} c_{9} + \frac{1}{383201280} c_{2} c_{6} c_{9}                                               \\
   & - \frac{1}{766402560} c_{3} c_{4} c_{10} - \frac{17}{4598415360} c_{2} c_{5} c_{10} - \frac{1}{4598415360} c_{3}^{2} c_{11} + \frac{1}{229920768} c_{2} c_{4} c_{11} + \frac{1}{191600640} c_{2} c_{3} c_{12}                                          \\
   & - \frac{1}{153280512} c_{2}^{2} c_{13} - \frac{1}{510935040} c_{8} c_{9} + \frac{1}{164229120} c_{7} c_{10} - \frac{1}{109486080} c_{6} c_{11} + \frac{1}{510935040} c_{5} c_{12}                                                                      \\
   & + \frac{19}{2299207680} c_{4} c_{13} - \frac{29}{1532805120} c_{3} c_{14} + \frac{1}{72990720} c_{2} c_{15} - \frac{41}{1532805120} c_{17} \bigg) \cdot v_5                                                                                            \\
   & - \int_M \bigg( \frac{1}{11412430848} c_{2}^{4} c_{3} c_{6} + \frac{1}{11412430848} c_{2}^{5} c_{7} + \frac{1}{11412430848} c_{2} c_{3}^{3} c_{6} + \frac{1}{3804143616} c_{2}^{2} c_{3} c_{4} c_{6}                                                   \\
   & + \frac{1}{11412430848} c_{2}^{3} c_{5} c_{6} - \frac{1}{11412430848} c_{2}^{2} c_{3}^{2} c_{7} - \frac{5}{11412430848} c_{2}^{3} c_{4} c_{7} + \frac{1}{5706215424} c_{2}^{3} c_{3} c_{8}                                                             \\
   & - \frac{1}{11412430848} c_{2}^{4} c_{9} - \frac{1}{11412430848} c_{3} c_{4}^{2} c_{6} - \frac{1}{11412430848} c_{3}^{2} c_{5} c_{6} - \frac{1}{5706215424} c_{2} c_{4} c_{5} c_{6}                                                                     \\
   & - \frac{1}{5706215424} c_{2} c_{3} c_{6}^{2} + \frac{5}{11412430848} c_{3}^{2} c_{4} c_{7} + \frac{5}{11412430848} c_{2} c_{4}^{2} c_{7} - \frac{1}{3804143616} c_{2} c_{3} c_{5} c_{7}                                                                \\
   & + \frac{1}{2853107712} c_{2}^{2} c_{6} c_{7} - \frac{5}{11412430848} c_{3}^{3} c_{8} - \frac{1}{2853107712} c_{2} c_{3} c_{4} c_{8} - \frac{1}{5706215424} c_{2}^{2} c_{5} c_{8}                                                                       \\
   & + \frac{1}{2853107712} c_{2} c_{3}^{2} c_{9} + \frac{1}{2853107712} c_{2}^{2} c_{4} c_{9} - \frac{1}{3804143616} c_{2}^{2} c_{3} c_{10} + \frac{1}{11412430848} c_{2}^{3} c_{11} + \frac{1}{11412430848} c_{5} c_{6}^{2}                               \\
   & + \frac{1}{2853107712} c_{5}^{2} c_{7} - \frac{1}{2853107712} c_{4} c_{6} c_{7} - \frac{5}{11412430848} c_{3} c_{7}^{2} - \frac{11}{11412430848} c_{4} c_{5} c_{8} + \frac{1}{713276928} c_{3} c_{6} c_{8}                                             \\
   & + \frac{5}{5706215424} c_{2} c_{7} c_{8} + \frac{11}{11412430848} c_{4}^{2} c_{9} - \frac{1}{2853107712} c_{3} c_{5} c_{9} - \frac{29}{11412430848} c_{2} c_{6} c_{9} - \frac{5}{5706215424} c_{3} c_{4} c_{10}                                        \\
   & + \frac{29}{11412430848} c_{2} c_{5} c_{10} + \frac{5}{11412430848} c_{3}^{2} c_{11} - \frac{1}{713276928} c_{2} c_{4} c_{11} - \frac{11}{5706215424} c_{2} c_{3} c_{12} + \frac{25}{11412430848} c_{2}^{2} c_{13}                                     \\
   & + \frac{5}{3804143616} c_{8} c_{9} - \frac{37}{11412430848} c_{7} c_{10} + \frac{41}{11412430848} c_{6} c_{11} - \frac{1}{2853107712} c_{5} c_{12} - \frac{37}{11412430848} c_{4} c_{13}                                                               \\
   & + \frac{47}{11412430848} c_{3} c_{14} - \frac{25}{11412430848} c_{2} c_{15} + \frac{37}{11412430848} c_{17} \bigg) \cdot v_6                                                                                                                           \\
   & -\int_M \bigg( \frac{1}{6469632000} c_{2}^{4} c_{3} c_{6} + \frac{1}{6469632000} c_{2}^{5} c_{7} + \frac{1}{6469632000} c_{2} c_{3}^{3} c_{6} + \frac{1}{2156544000} c_{2}^{2} c_{3} c_{4} c_{6} + \frac{1}{6469632000} c_{2}^{3} c_{5} c_{6}          \\
   & - \frac{1}{6469632000} c_{2}^{2} c_{3}^{2} c_{7} - \frac{1}{1293926400} c_{2}^{3} c_{4} c_{7} + \frac{1}{3234816000} c_{2}^{3} c_{3} c_{8} - \frac{1}{6469632000} c_{2}^{4} c_{9} - \frac{1}{6469632000} c_{3} c_{4}^{2} c_{6}                         \\
   & - \frac{1}{6469632000} c_{3}^{2} c_{5} c_{6} - \frac{1}{3234816000} c_{2} c_{4} c_{5} c_{6} - \frac{1}{3234816000} c_{2} c_{3} c_{6}^{2} - \frac{1}{4313088000} c_{3}^{2} c_{4} c_{7} + \frac{1}{1293926400} c_{2} c_{4}^{2} c_{7}                     \\
   & + \frac{7}{12939264000} c_{2} c_{3} c_{5} c_{7} + \frac{1}{1617408000} c_{2}^{2} c_{6} c_{7} + \frac{1}{4313088000} c_{3}^{3} c_{8} - \frac{1}{1617408000} c_{2} c_{3} c_{4} c_{8} - \frac{17}{12939264000} c_{2}^{2} c_{5} c_{8}                      \\
   & - \frac{1}{2587852800} c_{2} c_{3}^{2} c_{9} + \frac{7}{4313088000} c_{2}^{2} c_{4} c_{9} + \frac{1}{646963200} c_{2}^{2} c_{3} c_{10} - \frac{1}{539136000} c_{2}^{3} c_{11} + \frac{1}{6469632000} c_{5} c_{6}^{2}                                   \\
   & - \frac{1}{2587852800} c_{5}^{2} c_{7} - \frac{1}{1617408000} c_{4} c_{6} c_{7} + \frac{1}{4313088000} c_{3} c_{7}^{2} + \frac{1}{431308800} c_{4} c_{5} c_{8} - \frac{7}{12939264000} c_{3} c_{6} c_{8}                                               \\
   & + \frac{7}{12939264000} c_{2} c_{7} c_{8} - \frac{1}{431308800} c_{4}^{2} c_{9} + \frac{1}{2587852800} c_{3} c_{5} c_{9} - \frac{19}{12939264000} c_{2} c_{6} c_{9} - \frac{7}{12939264000} c_{3} c_{4} c_{10}                                         \\
   & - \frac{7}{12939264000} c_{2} c_{5} c_{10} - \frac{1}{4313088000} c_{3}^{2} c_{11} + \frac{59}{12939264000} c_{2} c_{4} c_{11} - \frac{11}{3234816000} c_{2} c_{3} c_{12} + \frac{1}{539136000} c_{2}^{2} c_{13}                                       \\
   & + \frac{1}{431308800} c_{8} c_{9} - \frac{1}{269568000} c_{7} c_{10} + \frac{1}{3234816000} c_{6} c_{11} + \frac{7}{1293926400} c_{5} c_{12} - \frac{1}{129392640} c_{4} c_{13}                                                                        \\
   & - \frac{23}{12939264000} c_{3} c_{14} + \frac{67}{12939264000} c_{2} c_{15} - \frac{41}{6469632000} c_{17} \bigg) \cdot v_7                                                                                                                            \\
   & - \int_M \bigg( \frac{1}{11496038400} c_{2}^{5} c_{3} c_{4} + \frac{1}{11496038400} c_{2}^{6} c_{5} + \frac{1}{5748019200} c_{2}^{2} c_{3}^{3} c_{4} + \frac{1}{2874009600} c_{2}^{3} c_{3} c_{4}^{2} - \frac{1}{5748019200} c_{2}^{3} c_{3}^{2} c_{5} \\
   & - \frac{1}{2299207680} c_{2}^{4} c_{4} c_{5} + \frac{1}{5748019200} c_{2}^{4} c_{3} c_{6} - \frac{1}{11496038400} c_{2}^{5} c_{7} - \frac{1}{11496038400} c_{3}^{3} c_{4}^{2} - \frac{1}{3832012800} c_{2} c_{3} c_{4}^{3}                             \\
   & + \frac{1}{22992076800} c_{2} c_{3}^{2} c_{4} c_{5} + \frac{1}{1916006400} c_{2}^{2} c_{4}^{2} c_{5} + \frac{1}{4598415360} c_{2}^{2} c_{3} c_{5}^{2} - \frac{1}{4598415360} c_{2} c_{3}^{3} c_{6} - \frac{1}{1277337600} c_{2}^{2} c_{3} c_{4} c_{6}  \\
   & - \frac{1}{22992076800} c_{2}^{3} c_{5} c_{6} + \frac{1}{22992076800} c_{2}^{2} c_{3}^{2} c_{7} + \frac{17}{22992076800} c_{2}^{3} c_{4} c_{7} + \frac{1}{7664025600} c_{2}^{3} c_{3} c_{8} - \frac{1}{3284582400} c_{2}^{4} c_{9}                     \\
   & - \frac{1}{11496038400} c_{4}^{3} c_{5} - \frac{1}{11496038400} c_{3} c_{4} c_{5}^{2} - \frac{1}{22992076800} c_{2} c_{5}^{3} + \frac{1}{2874009600} c_{3} c_{4}^{2} c_{6} + \frac{1}{4598415360} c_{3}^{2} c_{5} c_{6}                                \\
   & + \frac{1}{5748019200} c_{2} c_{4} c_{5} c_{6} + \frac{1}{2299207680} c_{2} c_{3} c_{6}^{2} - \frac{1}{11496038400} c_{3}^{2} c_{4} c_{7} - \frac{1}{1045094400} c_{2} c_{4}^{2} c_{7} - \frac{1}{4598415360} c_{2} c_{3} c_{5} c_{7}                  \\
   & - \frac{17}{22992076800} c_{2}^{2} c_{6} c_{7} - \frac{1}{5748019200} c_{2} c_{3} c_{4} c_{8} + \frac{1}{1437004800} c_{2}^{2} c_{5} c_{8} + \frac{1}{4598415360} c_{2} c_{3}^{2} c_{9} + \frac{1}{2299207680} c_{2}^{2} c_{4} c_{9}                   \\
   & - \frac{1}{1094860800} c_{2}^{2} c_{3} c_{10} + \frac{1}{1277337600} c_{2}^{3} c_{11} - \frac{1}{4598415360} c_{5} c_{6}^{2} + \frac{1}{22992076800} c_{5}^{2} c_{7} + \frac{1}{1149603840} c_{4} c_{6} c_{7}                                          \\
   & + \frac{1}{5748019200} c_{3} c_{7}^{2} - \frac{1}{1642291200} c_{4} c_{5} c_{8} - \frac{1}{4598415360} c_{3} c_{6} c_{8} - \frac{1}{1149603840} c_{2} c_{7} c_{8} + \frac{1}{11496038400} c_{4}^{2} c_{9}                                              \\
   & - \frac{1}{4598415360} c_{3} c_{5} c_{9} + \frac{17}{11496038400} c_{2} c_{6} c_{9} + \frac{23}{22992076800} c_{3} c_{4} c_{10} + \frac{1}{11496038400} c_{2} c_{5} c_{10}                                                                             \\
   & - \frac{29}{11496038400} c_{2} c_{4} c_{11} + \frac{17}{22992076800} c_{2} c_{3} c_{12} + \frac{1}{5748019200} c_{2}^{2} c_{13} - \frac{13}{7664025600} c_{8} c_{9} + \frac{59}{22992076800} c_{7} c_{10}                                              \\
   & + \frac{1}{4598415360} c_{6} c_{11} - \frac{61}{22992076800} c_{5} c_{12} + \frac{53}{22992076800} c_{4} c_{13} + \frac{1}{522547200} c_{3} c_{14} - \frac{61}{22992076800} c_{2} c_{15}                                                               \\
   & + \frac{61}{22992076800} c_{17} \bigg) \cdot v_8.
\end{align*}

\section{Chern numbers of $Y_1$ and $Y_2$} \label{app:chernnum}
\tablefirsthead{\Hline{1pt}
  & $Y_1$                 & $Y_2$           \rule[-8pt]{0pt}{20pt}      \\
  \hline}
\tablehead{\Hline{1pt}
  & $Y_1$                 & $Y_2$              \rule[-8pt]{0pt}{20pt}    \\
  \hline}
\tabletail{\Hline{1pt}}
\tablelasttail{\Hline{1pt}}
\begin{center}
  \begin{xtabular}{ccc}
    $\int \chern_{17}$                                                & $ -12566964323536824$ & $ -12566964323536824$ \\
    \hline
    $\int \chern_{15} \chern_{ 2}$                                    & $ -16120276856522016$ & $ -17771184345938256$ \\
    \hline
    $\int \chern_{14} \chern_{ 3}$                                    & $ -10418778246653592$ & $ -12069685736069832$ \\
    \hline
    $\int \chern_{13} \chern_{ 4}$                                    & $ -16955718164205288$ & $ -20118448416029328$ \\
    \hline
    $\int \chern_{13} \chern_{ 2} \chern_{ 2}$                        & $ -20017380823793664$ & $ -23731960802953824$ \\
    \hline
    $\int \chern_{12} \chern_{ 5}$                                    & $ -17629884014664912$ & $ -21561264608131872$ \\
    \hline
    $\int \chern_{12} \chern_{ 3} \chern_{ 2}$                        & $ -12704878750253568$ & $ -15676286308648848$ \\
    \hline
    $\int \chern_{11} \chern_{ 6}$                                    & $ -19868302036394160$ & $ -24734418956664120$ \\
    \hline
    $\int \chern_{11} \chern_{ 4} \chern_{ 2}$                        & $ -20280725863792512$ & $ -25380656468947872$ \\
    \hline
    $\int \chern_{11} \chern_{ 3} \chern_{ 3}$                        & $ -7981790182965696$  & $ -10210025083273296$ \\
    \hline
    $\int \chern_{11} \chern_{ 2} \chern_{ 2} \chern_{ 2}$            & $ -23947738758721536$ & $ -29939834074648176$ \\
    \hline
    $\int \chern_{10} \chern_{ 7}$                                    & $ -20981648718607848$ & $ -26473467597046848$ \\
    \hline
    $\int \chern_{10} \chern_{ 5} \chern_{ 2}$                        & $ -20643421387999488$ & $ -26407977811451568$ \\
    \hline
    $\int \chern_{10} \chern_{ 4} \chern_{ 3}$                        & $ -12600241652331024$ & $ -16282668831769224$ \\
    \hline
    $\int \chern_{10} \chern_{ 3} \chern_{ 2} \chern_{ 2}$            & $ -14878855623247872$ & $ -19207843059552912$ \\
    \hline
    $\int \chern_{9} \chern_{ 8}$                                     & $ -21640311771607104$ & $ -27477296908130904$ \\
    \hline
    $\int \chern_{9} \chern_{ 6} \chern_{ 2}$                         & $ -22725369462356352$ & $ -29363935511099952$ \\
    \hline
    $\int \chern_{9} \chern_{ 5} \chern_{ 3}$                         & $ -12667491348925152$ & $ -16674977592147072$ \\
    \hline
    $\int \chern_{9} \chern_{ 4} \chern_{ 4}$                         & $ -19648034664639552$ & $ -25550516571066072$ \\
    \hline
    $\int \chern_{9} \chern_{ 4} \chern_{ 2} \chern_{ 2}$             & $ -23203089358654464$ & $ -30140610675346224$ \\
    \hline
    $\int \chern_{9} \chern_{ 3} \chern_{ 3} \chern_{ 2}$             & $ -9130259521529856$  & $ -12130917989089536$ \\
    \hline
    $\int \chern_{9} \chern_{ 2} \chern_{ 2} \chern_{ 2} \chern_{ 2}$ & $ -27402915124445184$ & $ -35555242134916224$ \\
    \hline
    $\int \chern_{8} \chern_{ 7} \chern_{ 2}$                         & $ -23382771727096320$ & $ -30406320596710800$ \\
    \hline
    $\int \chern_{8} \chern_{ 6} \chern_{ 3}$                         & $ -13755858825604608$ & $ -18229439098901448$ \\
    \hline
    $\int \chern_{8} \chern_{ 5} \chern_{ 4}$                         & $ -19485284047106400$ & $ -25723444707610560$ \\
    \hline
    $\int \chern_{8} \chern_{ 5} \chern_{ 2} \chern_{ 2}$             & $ -23011236776890368$ & $ -30344762098402848$ \\
    \hline
    $\int \chern_{8} \chern_{ 4} \chern_{ 3} \chern_{ 2}$             & $ -14045147931715392$ & $ -18714726521158752$ \\
    \hline
    $\int \chern_{8} \chern_{ 3} \chern_{ 3} \chern_{ 3}$             & $ -5526623743186752$  & $ -7534042519741392$  \\
    \hline
    $\int \chern_{8} \chern_{ 3} \chern_{ 2} \chern_{ 2} \chern_{ 2}$ & $ -16587436762066944$ & $ -22077022229648784$ \\
    \hline
    $\int \chern_{7} \chern_{ 7} \chern_{ 3}$                         & $ -13911339130144848$ & $ -18499829570435568$ \\
    \hline
    $\int \chern_{7} \chern_{ 6} \chern_{ 4}$                         & $ -20797455706552584$ & $ -27558359677576944$ \\
    \hline
    $\int \chern_{7} \chern_{ 6} \chern_{ 2} \chern_{ 2}$             & $ -24561353964006144$ & $ -32509349693941824$ \\
    \hline
    $\int \chern_{7} \chern_{ 5} \chern_{ 5}$                         & $ -18992797533920064$ & $ -25380945525785664$ \\
    \hline
    $\int \chern_{7} \chern_{ 5} \chern_{ 3} \chern_{ 2}$             & $ -13690228978704384$ & $ -18466904154925584$ \\

    \hline
    $\int \chern_{7} \chern_{ 4} \chern_{ 4} \chern_{ 2}$                                                 & $ -21236959469326080$ & $ -28292854220232480$ \\
    \hline
    $\int \chern_{7} \chern_{ 4} \chern_{ 3} \chern_{ 3}$                                                 & $ -8355975668453328$  & $ -11390667286338288$ \\
    \hline
    $\int \chern_{7} \chern_{ 4} \chern_{ 2} \chern_{ 2} \chern_{ 2}$                                     & $ -25082627990780928$ & $ -33375875438854368$ \\
    \hline
    $\int \chern_{7} \chern_{ 3} \chern_{ 3} \chern_{ 2} \chern_{ 2}$                                     & $ -9868440929252352$  & $ -13437272003300352$ \\
    \hline
    $\int \chern_{7} \chern_{ 2} \chern_{ 2} \chern_{ 2} \chern_{ 2} \chern_{ 2}$                         & $ -29626288313401344$ & $ -39372018689971584$ \\
    \hline
    $\int \chern_{6} \chern_{ 6} \chern_{ 5}$                                                             & $ -20045163757887552$ & $ -26794314133031232$ \\
    \hline
    $\int \chern_{6} \chern_{ 6} \chern_{ 3} \chern_{ 2}$                                                 & $ -14448969306563136$ & $ -19495532608483536$ \\
    \hline
    $\int \chern_{6} \chern_{ 5} \chern_{ 4} \chern_{ 2}$                                                 & $ -20468820683488704$ & $ -27509835217218624$ \\
    \hline
    $\int \chern_{6} \chern_{ 5} \chern_{ 3} \chern_{ 3}$                                                 & $ -8053707546698976$  & $ -11076233891550336$ \\
    \hline
    $\int \chern_{6} \chern_{ 5} \chern_{ 2} \chern_{ 2} \chern_{ 2}$                                     & $ -24175413086373888$ & $ -32452298318020608$ \\
    \hline
    $\int \chern_{6} \chern_{ 4} \chern_{ 4} \chern_{ 3}$                                                 & $ -12493244792291808$ & $ -16968762117951048$ \\
    \hline
    $\int \chern_{6} \chern_{ 4} \chern_{ 3} \chern_{ 2} \chern_{ 2}$                                     & $ -14755434031446528$ & $ -20017535250734928$ \\
    \hline
    $\int \chern_{6} \chern_{ 3} \chern_{ 3} \chern_{ 3} \chern_{ 2}$                                     & $ -5805373830183936$  & $ -8060435949465216$  \\
    \hline
    $\int \chern_{6} \chern_{ 3} \chern_{ 2} \chern_{ 2} \chern_{ 2} \chern_{ 2}$                         & $ -17428190600921088$ & $ -23614029246446208$ \\
    \hline
    $\int \chern_{5} \chern_{ 5} \chern_{ 5} \chern_{ 2}$                                                 & $ -18692288718173184$ & $ -25337872758579264$ \\
    \hline
    $\int \chern_{5} \chern_{ 5} \chern_{ 4} \chern_{ 3}$                                                 & $ -11408947399354176$ & $ -15629543395106976$ \\
    \hline
    $\int \chern_{5} \chern_{ 5} \chern_{ 3} \chern_{ 2} \chern_{ 2}$                                     & $ -13474616955015168$ & $ -18437767402667328$ \\
    \hline
    $\int \chern_{5} \chern_{ 4} \chern_{ 4} \chern_{ 4}$                                                 & $ -17698706221829376$ & $ -23944319006461056$ \\
    \hline
    $\int \chern_{5} \chern_{ 4} \chern_{ 4} \chern_{ 2} \chern_{ 2}$                                     & $ -20904367685689344$ & $ -28246298751476064$ \\
    \hline
    $\int \chern_{5} \chern_{ 4} \chern_{ 3} \chern_{ 3} \chern_{ 2}$                                     & $ -8224242734108160$  & $ -11373973294302720$ \\
    \hline
    $\int \chern_{5} \chern_{ 4} \chern_{ 2} \chern_{ 2} \chern_{ 2} \chern_{ 2}$                         & $ -24691976897396736$ & $ -33321125005612416$ \\
    \hline
    $\int \chern_{5} \chern_{ 3} \chern_{ 3} \chern_{ 3} \chern_{ 3}$                                     & $ -3235804304358912$  & $ -4580559535500672$  \\
    \hline
    $\int \chern_{5} \chern_{ 3} \chern_{ 3} \chern_{ 2} \chern_{ 2} \chern_{ 2}$                         & $ -9713731756621824$  & $ -13417639482731904$ \\
    \hline
    $\int \chern_{5} \chern_{ 2} \chern_{ 2} \chern_{ 2} \chern_{ 2} \chern_{ 2} \chern_{ 2}$             & $ -29167396545626112$ & $ -39307625506351872$ \\
    \hline
    $\int \chern_{4} \chern_{ 4} \chern_{ 4} \chern_{ 3} \chern_{ 2}$                                     & $ -12758888683201728$ & $ -17424186820264368$ \\
    \hline
    $\int \chern_{4} \chern_{ 4} \chern_{ 3} \chern_{ 3} \chern_{ 3}$                                     & $ -5019681599497920$  & $ -7016824292191920$  \\
    \hline
    $\int \chern_{4} \chern_{ 4} \chern_{ 3} \chern_{ 2} \chern_{ 2} \chern_{ 2}$                         & $ -15070492635303936$ & $ -20554809649453296$ \\
    \hline
    $\int \chern_{4} \chern_{ 3} \chern_{ 3} \chern_{ 3} \chern_{ 2} \chern_{ 2}$                         & $ -5928734983747584$  & $ -8277654181388544$  \\
    \hline
    $\int \chern_{4} \chern_{ 3} \chern_{ 2} \chern_{ 2} \chern_{ 2} \chern_{ 2} \chern_{ 2}$             & $ -17801848444551168$ & $ -24247859797369728$ \\
    \hline
    $\int \chern_{3} \chern_{ 3} \chern_{ 3} \chern_{ 3} \chern_{ 3} \chern_{ 2}$                         & $ -2332510279839744$  & $ -3333926642649984$  \\
    \hline
    $\int \chern_{3} \chern_{ 3} \chern_{ 3} \chern_{ 2} \chern_{ 2} \chern_{ 2} \chern_{ 2}$             & $ -7002794034462720$  & $ -9765015249070080$  \\
    \hline
    $\int \chern_{3} \chern_{ 2} \chern_{ 2} \chern_{ 2} \chern_{ 2} \chern_{ 2} \chern_{ 2} \chern_{ 2}$ & $ -21029331652313088$ & $ -28604369215531008$ \\
  \end{xtabular}
\end{center}

\bibliography{references}
\bibliographystyle{amsalpha}
\end{document}